\documentclass[a4paper,12pt]{article}
\usepackage{amssymb}
\usepackage[T2A]{fontenc}
\usepackage[cp1251]{inputenc}
\usepackage{amsmath}

\usepackage[dvips]{graphicx}
\usepackage[english,russian]{babel}

\newtheorem{theorem}{Теорема}

\newtheorem{definitionhead}[theorem]{Определение}

\author{А.~Я.~Белов (МФТИ), Г.~О.~Шнайдер (ДНТТМ)}

\title{Об адаптированном курсе математики для кружковцев-химиков}

\begin{document}

\maketitle

  \begin{abstract}

    Статья посвящена опыту преподавания математики кружковцам-химикам. Перекликается с опытом И.С.Рубанова и С.Е.Рукшина. В условиях (двух)семестрового кружка ничему нельзя научить, но можно показать силу красоты и как математическая красота связана с физическом смыслом на одном примере. А далее человек сам ориентируется.

\end{abstract}

\section{Введение}

     Важность  развития  математической  культуры  у  специалистов
самых  разных  отраслей  несомненна.  Однако   современная   школа
зачастую прививает  отвращение  к  математике.  В  школьном  курсе
внимание акцентируется на ошибках,  том,  что  не  следует  делать
(примером  чему  служат  бесконечные  ``ОДЗ''),   а   само   понятие
математической  строгости  фетишируется.  Все  это   не   вызывает
энтузиазма у  учащихся  и  мы  получаем  извращенное  отношение  к
предмету - как к своду формальных правил, применяемых механически.
Именно такое  отношение  прежде  всего  и  вызывает  разговоры  об
``отсутствии математических способностей'', в то время как  практика
показывает,  что  особенностью   многих   (настоящих)   отличников
является не интеллект, а  простой  здравый  смысл  в  отношении  к
предмету. Об этом подробнее писалось в работе \cite{BelovYavich}, данная статья с ней перекликается.
Нам созвучны мысли, изложенные в работах \cite{RukshinMaximov,RukshinMaximov1,RukshinMaximov2,RukshinMaximov3}.
О преодолении возникающих
методических проблем в процессе преподавании индукции написано в замечательной статье
И.~С.~Рубанова \cite{Rubanov}.

     Потребность в преодолении этого отношения и вызвало  к  жизни
адаптивный курс математики, разработанный в рамках  сотрудничества
лабораторий математики  и  коллоидной  химии.

\section{Адаптивный курс математики}

Основная  концепция
курса состоит в том,  что  хотя  и  невозможно  дать  сколь-нибудь
систематические   знания   в   рамках   семестрового   или    даже
двухсеместрового курса, но можно преодолеть школьное  отношение  к
предмету и выработать правильное, повернуть сознание учащегося. Таким образом, целями и задачами
курса являются:

\begin{itemize}
  \item преодоление  формализма  в  отношении   к   математике   и
привлечение здравого смысла учащегося к предмету.

  \item показ связи между  математическим  и  физическим  смыслами.
Ученик должен понимать, что  написав  математическую  формулу,  он
сделал некоторую физическую гипотезу, и надо понять, какую именно.
  \item развитие эстетического чувства учащихся и общей культуры.
  \item сообщение учащимся знаний, поддерживающих  курс  коллоидной
химии. (теория вероятностей, понятие  интеграла  и  производной  и
т.д.)
  \item приучить  учащихся  к  физическому  уровню   строгости   в
математических рассуждениях.
\end{itemize}

 \subsection{Основные средства решения этих задач}

\begin{itemize}
  \item через  разбор  примеров,  имеющих  несомненное  прикладное
значение для учащегося демонстрировать связь  между  физическим  и
математическим смыслами.
  \item исторический    рассказ    об    изобретении    некоторых
математических  понятий  ``практиками''.  При   этом   акцентируется
внимание   на   физическом   уровне   строгости   рассуждений   и
противопоставлении его  школьному  стилю.  (Например,  рассказ  об
изобретении логарифмов.)
  \item рассказ    эстетически    значимых    сюжетов,    имеющих
общекультурное и общефилосовское  значение  (например,  логические
парадоксы).
\end{itemize}

\subsection{Список  основных  тем,  которые  читались  в
разное время}

\begin{itemize}
  \item Логические парадоксы (2 часа).
  \item Понятия {\it скорости} и {\it производной} (4 часа).
  \item Процедура интегрирования (4 часа).
  \item Понятие {\it плотности} (2 часа).
  \item История создания логарифмов (4 часа).
  \item Техника счета пределов (4 часа).
  \item Техника дифференцирования (4 часа).
  \item Комбинаторика (4 часа).
  \item Понятие {\it вероятности} (2 часа).
  \item Биномиальное распределение (2 часа).
  \item Геометрические вероятности (2 часа).
  \item Пуассоновское распределение (2 часа).
  \item Случайные блуждания. Матожидание и дисперсия (4 часа).
  \item Оценка скорости диффузии (2 часа).
  \item О Гауссовом распределении (2 часа).
  \item Распределение Максвелла (2 часа).
  \item Оценка скорости химических реакций (2 часа).
\end{itemize}

\section{Заключение}

     Сотрудничество лабораторий ДНТТМ продолжалось более 8 лет, велись совместные кружки, на основании
     опыта которых и возникла данная статья (Г.~О.~Шнайдер долгие годы работал заведующим лаборатории коллоидной химии).
     За  время  работы  были   достигнуты   определенные   успехи.
Некоторые  кружковцы  --   химики   заняли   призовые   места   на
математических  и  физических  олимпиадах  (районная,   городская,
международный турнир  городов),  математические  доклады  учащихся
были отмечены премиями на конференциях ``Поиск'' и на  международной
конференции памяти Чижевского (г.Обнинск). (Доклады учащихся были посвящены оценкам скорости химических реакций и распределению Максвелла. См., например, \cite{Kuligin})
Отметим,  что  успех  на
олимпиаде или на конференции сам по себе,  никогда  не  был  нашей
основной целью.  Это  лишь  следствие  изменившегося  отношения  к
предмету. В дальнейшем предполагалось продолжить работу  в  рамках
совместного  проекта  трех  лабораторий  (математики,   коллоидной
химии,   вычислительной   техники).   Проект   посвящен   проблеме
разделения    эмульсии    и    может    иметь     непосредственное
народохозяйственное значение.

\section{Приложение. Оценки скорости химических реакций и распределение Максвелла-Больцмана}

{\bf Введение.} Давайте удивимся. Когда температура у человека 38 градусов он чувствует себя совсем по-другому чем когда 37. Но в абсолютном выражении это 311 или 310 градусов Кельвина, т.е. средняя энергия молекулы отличается на $0.3\%$. Откуда такое значение такой незначительной величины? Чтобы это понять, надо разобраться с распределением молекул по скоростям. Химические реакции обеспечиваются быстрыми частицами и оказывается, количество быстрых частиц изменится заметно, но как именно -- это тема данного курса.

{\bf 1.} Чтобы разобраться, как распределены молекулы по скоростям, нужна концепция {\it плотности вероятности}, ведь вероятность того, что скорость именно такая равна нулю!
Обсуждаем концепции {\it плотности},  а также {\it плотности вероятности}.

{\bf 2.} Постулаты Максвелла: 1) Плотность вероятности не зависит от направления вектора скорости (изотропность пространства). 2) Компоненты скоростей $(V_x,V_y,V_z)$ распределены независимо. Обсуждается подробно понятие {\it независимости событий} для классических вероятностей и для плотностей.

Формулируется итог: $\varphi(\vec{V}^2)=F(\vec{V})=\phi(V_x^2)\phi(V_y^2)\phi(V_z^2)=\varphi(V_x^2+V_y^2+V_z^2)$. Обозначим $a=V_x^2$, $b=v_y^2$, $c=v_z^2$. Имеем:
$$\varphi(a+b+c)=\phi(a)\phi(b)\phi(c), \quad \ln(\varphi(a+b+c))=\ln(\phi(a))+\ln(\phi(b))+\ln(\phi(c))$$

при любых $a,b,c$. Решаем функциональное уравнение. Попутно обсуждаем понятие {\it функции} как {\it закона соответствия}, который нужно найти (исторически понятие  {\it функции} возникло, когда потребовалось решить функциональное уравнение, т.е. из множества законов соответствия выбрать один правильный.

В итоге имеем: $F(\vec{V})=C_1\exp(C_2\vec{V}^2)$.

{\bf 3.} \ Находим константы $C_1, C_2$. $\int_{\vec{V}} F(\vec{V})=1$, т.е. $$\int_{-\infty}^{\infty}\int_{-\infty}^{\infty}\int_{-\infty}^{\infty}C_2\exp(C_1V_x^2)dV_x\exp(C_1V_y^2)dV_y\exp(C_1V_z^2)dV_z=1.$$

Понятие {\it интеграла}. Как известно, Байкал создан для нужд целлюлозно-бумажного комбината. Следует определить запасы престной воды в озере. Пусть оно замерзло. Разобем поверхность озера на прямоугольнички со сторонвми $dx$ и $dy$. Буровая установка делает прямоугольное отверстие, и льда достается $h(x,y)dxdy$ из каждого столбика. Общий запас воды равен $\sum_{(x,y)}h(x,y)dxdy$. При удлинении знака $S$ (Summa onium) получается интеграл. Пересыхает Аральское море. уровень уменьшился на $dh$, а площадь поверхности $S(h)$. Объем высохшей воды равен $S(h)dh$, а общий объем $\int_h S(h)dh=\int\int h(x,y)dxdy$. Интеграл от плотности.
Вопросы сходимости игнорируются ради наглядности. Интеграл от плотности есть масса. Приводятся примеры.

{\bf Интеграл Гаусса.} $J=\int_{-\infty}^\infty e^{-x^2}dx=\sqrt{\pi}$. Тогда $$J^2=\int_{-\infty}^\infty e^{-x^2}dx\int_{-\infty}^\infty e^{-y^2}dy=\int_{-\infty}^\infty e^{-(x^2+y^2)}dxdy$$ Двумя способами вычисляется площадь под поверхностью $e^{-(x^2+y^2)}=e^{-r^2}$. С одной стороны это $J^2$, (разбиваем на слои вдоль оси $(0X)$, с другой стороны это $\int 2\pi re^{-r^2}dr$  в силу разбиения на цилиндрические слои.

В итоге имеем $J^2=\pi$ или $J=\sqrt{\pi}$. Доказываем что $C_1=-C<0$, $C_2=(C/\pi)^{3/2}$ из равенства единице интеграла от плотности вероятности. Итак, $$F(\vec{V})=(C/\pi)^{3/2}e^{-C\vec{V}^2}.$$

{\bf 3.} Определяем константу $C$ исходя из средней энергии частицы, которая равна {\it температуре}. Энергия частицы со скоростью $\vec{V}$ равна $m\vec{V}^2/2$ а средняя энергия равна $$m\int_{\vec{V}}\vec{V}^2/2(C/\pi)^{3/2}e^{-C\vec{V}^2}dV_xdV_ydV_z=m\int_{\vec{V}}(V_x^2+V_y^2+V_z^2)/2(C/\pi)^{3/2}e^{-C\vec{V}^2}dV_xdV_ydV_z=$$
$$=3/2m(C/\pi)^{3/2}\int_{-\infty}^\infty V_x^2e^{-CV_x^2}dV_x\int_{-\infty}^\infty e^{-CV_y^2}dV_y
\int_{-\infty}^{\infty} e^{-CV_z^2}dV_z=$$
$$=m\sqrt{C/\pi}\int_{-\infty}^{+\infty}V_x^2e^{-CV_x^2}dV_x=\frac{3/2m}{C}=3/2kT$$

Откуда $C=m/kT$.

{\bf 4.} Использование свойств распределения. Время пробега молекулы -- порядка наносекунды, если нужной энергией обладают $10^{-9}$ молекул, то реакция проходит за время порядка секунды. Если же нужной энергией обладают $10^{-12}$ молекул (если $e^{-\lambda}=10^{-12}$, то $\lambda\sim 30$ то реакция проходит за время порядка 1000 секунд. В этом случае легко показать, что пр и увеличении средней энергии на $0.3\%$ количество частиц с нужной энергией возрастет весьма значительно, примерно на $10\%$.

{\bf 5.} По ходу дела могут обсуждаться случайные блуждания, оцениваться длина пробега в зависимости от диффузии и т.д. В зависимости от времени и силы школьников.

%\vspace*{-2\baselineskip}
\vspace*{-\baselineskip}
\enlargethispage{2\baselineskip}


\begin{thebibliography}{199}

\bibitem{BelovYavich}
{\sl Белов Алексей Яковлевич, Роман Явич}, {\it Проблемы одаренности и стадийность математического обучения}, Математическое образование, 3:1, Январь-Март (2010), 2--5, Engl. transl. {\sl A.Kanel-Belov, R. Yavich} (2012) {\it Staging of Mathematical Education},// Athens Institute for Education and Research// ATINER's Conference Paper Series MAT2012-0176// Athens, Greece, June 2012, pp 5--9 // ISSN 2241--2891, http://www.atiner.gr/papers/MAT2012-0176.pd

\bibitem{Kuligin} {\sl А.~К.~Кулыгин} (Рук. Белов А.~Я.)
{\it Распределение Максвелла.} Поиск -- 93, материалы конференции, физика, математика, астрономия, экология, 1993, Москва, ДНТТМ, стр. 23--29 


\bibitem{Rubanov}
{\sl Рубанов И.~С.} {\it Как обучать методу математической индукции?} Математика в школе. -- N1, 1996, стр. 14--20.

\bibitem{RukshinMaximov}
{\sl Рукшин С.~Е. Максимов Д.~В.} {\it
Структура интуиции у студентов-математиков и ее использование в решении задач.}
Некоторые актуальные
проблемы современной математики и математического образования. – СПб.:
БАН, 2005, с.152–-158.

\bibitem{RukshinMaximov1}
{\sl Рукшин С.~Е. Максимов Д.~В.} {\it
  Сравнительный анализ различных
аспектов процесса решения задач как
модели математического творчества -– $II$}
 РГПУ им. Герцена, некоторые актуальные проблемы современной математики и матем. образования,
 Матер. научн. конф. Герценовск. чтения,  $LIX$
16--21 апреля 2007, СПБ 2007, стр. 239--245


\bibitem{RukshinMaximov2}
{\sl Рукшин С.~Е. Максимов Д.~В.} {\it
  Сравнительный анализ различных
аспектов процесса решения задач как
модели математического творчества -– $I$}
РГПУ им. Герцена, некоторые актуальные проблемы современной математики и матем. образования,
 Матер. научн. конф. Герценовск. чтения,$LX$
17--22 апреля 2007, СПБ 2006,
 стр. 168--172


\bibitem{RukshinMaximov3}
{\sl Рукшин С.~Е. Максимов Д.~В.} {\it
  Сравнительный анализ различных
аспектов процесса решения задач как
модели математического творчества -– $III$}
РГПУ им. Герцена, некоторые актуальные проблемы современной математики и матем. образования,
 Матер. научн. конф. Герценовск. чтения, $LXII$
13--18 апреля 2009, СПБ 2009,
 стр. 196--199




\end{thebibliography}
\end{document}